\documentclass[reqno]{amsart}
\usepackage{amssymb}
\usepackage[pagebackref,hypertex]{hyperref}
\allowdisplaybreaks[4]
\theoremstyle{plain}
\newtheorem{thm}{Theorem}
\newtheorem{lem}{Lemma}

\theoremstyle{remark}
\newtheorem{rem}{Remark}
\DeclareMathOperator{\td}{d}
\date{Completed on 24 November 2008 in VU's Student Village}
\date{Revised from Open-TJM-2003.tex on 28 November 2010}
\date{}

\begin{document}

\title[Monotonicity of a function involving gamma functions]
{Complete monotonicity of a function involving a ratio of gamma functions and applications}

\author[F. Qi]{Feng Qi}
\address[F. Qi]{Department of Mathematics, College of Science, Tianjin Polytechnic University, Tianjin City, 300160, China; School of Mathematics and Informatics, Henan Polytechnic University, Jiaozuo City, Henan Province, 454010, China}
\email{\href{mailto: F. Qi <qifeng618@gmail.com>}{qifeng618@gmail.com}, \href{mailto: F. Qi <qifeng618@hotmail.com>}{qifeng618@hotmail.com}, \href{mailto: F. Qi <qifeng618@qq.com>}{qifeng618@qq.com}}
\urladdr{\url{http://qifeng618.wordpress.com}}

\author[B.-N. Guo]{Bai-Ni Guo}
\address[B.-N. Guo]{School of Mathematics and Informatics, Henan Polytechnic University, Jiaozuo City, Henan Province, 454010, China}
\email{\href{mailto: B.-N. Guo <bai.ni.guo@gmail.com>}{bai.ni.guo@gmail.com}, \href{mailto: B.-N. Guo <bai.ni.guo@hotmail.com>}{bai.ni.guo@hotmail.com}}

\begin{abstract}
In the paper, necessary and sufficient conditions are presented for a function involving a ratio of gamma functions to be logarithmically completely monotonic. This extends and generalizes the main result in [\emph{Inequalities and monotonicity for the ratio of gamma functions}, Taiwanese J. Math. \textbf{7} (2003), no.~2, 239\nobreakdash--247.] and others. As applications, several inequalities involving the volume of the unit ball in $\mathbb{R}^n$ are derived, which refine, generalize and extend some known inequalities.
\end{abstract}



\thanks{The first author was partially supported by the China Scholarship Council and the Science Foundation of Tianjin Polytechnic University}

\thanks{This paper was typeset using \AmS-\LaTeX}

\maketitle

\section{Introduction}

Recall from \cite{Atanassov, compmon2} that a positive real-valued function $f(x)$ is said to be logarithmically completely monotonic on an interval $I\subseteq\mathbb{R}$ if it has derivatives of all orders on $I$ and its logarithm $\ln f$ satisfies
$$
0\le(-1)^k[\ln f(x)]^{(k)}<\infty
$$
for $k\in\mathbb{N}$ on $I$. For more properties of this class of functions, please refer to~\cite{CBerg}.
\par
It is general knowledge that the classical Euler gamma function $\Gamma(x)$ may be defined for $x>0$ by
\begin{equation}\label{egamma}
\Gamma(x)=\int^\infty_0t^{x-1} e^{-t}\td t.
\end{equation}
The logarithmic derivative of $\Gamma(x)$, denoted by $\psi(x)=\frac{\Gamma'(x)}{\Gamma(x)}$, is called the psi or digamma function, and $\psi^{(k)}(x)$ for $k\in\mathbb{N}$ are called the polygamma functions. It is well known that these functions are fundamental and that they have much extensive applications in mathematical sciences.
\par
In~\cite[Theorem~2]{Guo-Qi-TJM-03.tex}, the following monotonicity was established: The function
\begin{equation}\label{fun-ori}
\frac{[{\Gamma(x+y+1)}/{\Gamma(y+1)}]^{1/x}}{x+y+1}
\end{equation}
is decreasing with respect to $x\ge1$ for fixed $y\ge0$.
Consequently, for positive real numbers $x\ge1$ and $y\ge0$, we have
\begin{equation}\label{new3}
\frac{x+y+1}{x+y+2}\le\frac{[\Gamma(x+y+1)/\Gamma(y+1)]^{1/x}}
{[\Gamma(x+y+2)/\Gamma(y+1)]^{1/(x+1)}}.
\end{equation}
\par
In~\cite{Extension-TJM-2003.tex}, the function~\eqref{fun-ori} was proved to be logarithmically completely monotonic with respect to $x\in(0,\infty)$ for $y\ge0$, and so is its reciprocal for $-1<y\le-\frac12$. Consequently, the inequality~\eqref{new3} is valid for $(x,y)\in(0,\infty)\times[0,\infty)$ and reversed for $(x,y)\in(0,\infty)\times\bigl(-1,-\frac12\bigr]$.
\par
For $(x,y)\in(0,\infty)\times[0,\infty)$ and $\alpha\in[0,\infty)$, the function
\begin{equation}\label{fun-F}
\frac{[{\Gamma(x+y+1)}/{\Gamma(y+1)}]^{1/x}}{(x+y+1)^\alpha}
\end{equation}
was proved in~\cite{Zhao-Chu-Jiang} to be strictly increasing (or decreasing, respectively) with respect to the single variable $x\in(0,\infty)$ if and only if $0\le\alpha\le\frac12$ (or $\alpha\ge1$, respectively), to be strictly increasing with respect to $y$ on $[0,\infty)$ if and only if $0\le\alpha\le1$, and to be logarithmically concave with respect to the $2$-variable $(x,y)\in(0,\infty)\times(0,\infty)$ if $0\le\alpha\le\frac14$.
\par
For given $y\in(-1,\infty)$ and $\alpha\in(-\infty,\infty)$, let
\begin{equation}\label{fun}
h_{\alpha,y}(x)=
\begin{cases}
\dfrac1{(x+y+1)^\alpha}\biggl[\dfrac{\Gamma(x+y+1)}{\Gamma(y+1)}\biggr]^{1/x}, &x\in(-y-1,\infty)\setminus\{0\};\\[0.8em]
\dfrac1{(y+1)^\alpha}\exp[\psi(y+1)],&x=0.
\end{cases}
\end{equation}
It is clear that the ranges of $x$, $y$ and $\alpha$ in the function $h_{\alpha,y}(x)$ extend the corresponding ones in the functions~\eqref{fun-ori} and~\eqref{fun-F} which were ever discussed in~\cite{Guo-Qi-TJM-03.tex, Extension-TJM-2003.tex, Zhao-Chu-Jiang}.
\par
The aim of this paper is to present necessary and sufficient conditions such that the function~\eqref{fun} or its reciprocal are logarithmically completely monotonic.
\par
Our main results may be stated as follows.

\begin{thm}\label{Open-TJM-2003-thm-1}
For $y>-1$, we have the following statements:
\begin{enumerate}
\item
the function~\eqref{fun} is logarithmically completely monotonic with respect to $x\in(-y-1,\infty)$ if and only if $\alpha\ge\max\bigl\{1,\frac1{y+1}\bigr\}$;
\item
if $\alpha\le\min\bigl\{1,\frac1{2(y+1)}\bigr\}$, the reciprocal of the function~\eqref{fun} is logarithmically completely monotonic with respect to $x\in(-y-1,\infty)$;
\item
the necessary condition for the reciprocal of the function~\eqref{fun} to be logarithmically completely monotonic with respect to $x\in(-y-1,\infty)$ is $\alpha\le1$.
\end{enumerate}
\end{thm}

As a ready consequence of monotonic results in Theorem~\ref{Open-TJM-2003-thm-1}, the following double inequality may be derived.

\begin{thm}\label{Open-TJM-2003-thm-ineq}
For $t>0$, $y+1>0$ and $x+y+1>0$, the double inequality
\begin{equation}\label{open-TJM-2003-ineq}
\biggl(\frac{x+y+1}{x+y+t+1}\biggr)^a <\frac{[\Gamma(x+y+1)/\Gamma(y+1)]^{1/x}}
{[\Gamma(x+y+t+1)/\Gamma(y+1)]^{1/(x+t)}} <\biggl(\frac{x+y+1}{x+y+t+1}\biggr)^b
\end{equation}
holds if $a\ge\max\bigl\{1,\frac1{y+1}\bigr\}$ and $b\le\min\bigl\{1,\frac1{2(y+1)}\bigr\}$.
\end{thm}

In order to show the applicability of Theorem~\ref{Open-TJM-2003-thm-ineq}, we derive the following double inequalities involving the $n$-dimensional volume
\begin{equation}
\Omega_n=\frac{\pi^{n/2}}{\Gamma(1+n/2)}
\end{equation}
of the unit ball $\mathbb{B}^n$ in $\mathbb{R}^n$.

\begin{thm}\label{vol-ball-ineq-thm}
For $n\in\mathbb{N}$, we have
\begin{gather}\label{sqrt-frac-n+2n+4}
\sqrt{\frac{n+2}{n+4}}\,<\frac{\Omega_{n+2}^{1/(n+2)}}{\Omega_n^{1/n}} <\sqrt[4]{\frac{n+2}{n+4}}\,,\\\label{vol-ball-1}
\frac1{\pi^{2/(n-2)n}}\sqrt{\frac{n+2}{n+4}}\, <\frac{\Omega_{n+2}^{1/n}}{\Omega_n^{1/(n-2)}}
<\frac1{\pi^{2/(n-2)n}}\sqrt[8]{\frac{n+2}{n+4}}\,,\\\label{ratio-Omega-n-n+1}
\sqrt{\frac{n+2}{n+3}}\, <\frac{\Omega_{n+1}^{1/(n+1)}}
{\Omega_{n}^{1/n}} <\sqrt[4]{\frac{n+2}{n+3}}\,.
\end{gather}
\end{thm}

In the final section, we will give several remarks about these three theorems.

\section{Lemmas}

In order to prove our main results, the following lemmas are needed.

\begin{lem}[{\cite[p.~107, Lemma~3]{theta-new-proof.tex-BKMS}}]\label{comp-thm-1}
For $x\in(0,\infty)$ and $k\in\mathbb{N}$, we have
\begin{equation}\label{qi-psi-ineq-1}
\ln x-\frac1x<\psi(x)<\ln x-\frac1{2x}
\end{equation}
and
\begin{equation}\label{qi-psi-ineq}
\frac{(k-1)!}{x^k}+\frac{k!}{2x^{k+1}}< (-1)^{k+1}\psi^{(k)}(x)<\frac{(k-1)!}{x^k}+\frac{k!}{x^{k+1}}.
\end{equation}
\end{lem}

\begin{lem}\label{qi-psi-ineq-beta-lem}
For $x\in(0,\infty)$ and $k\in\mathbb{N}$, we have
\begin{equation}\label{qi-psi-ineq-beta-1}
\ln\biggr(x+\frac12\biggl)-\frac1x<\psi(x)<\ln(x+1)-\frac1x
\end{equation}
and
\begin{equation}\label{qi-psi-ineq-beta-2}
\frac{(k-1)!}{(x+1)^k}+\frac{k!}{x^{k+1}}< (-1)^{k+1}\psi^{(k)}(x) <\frac{(k-1)!}{(x+1/2)^k}+\frac{k!}{x^{k+1}}.
\end{equation}
\end{lem}

\begin{proof}
In \cite[Theorem~1]{Guo-Qi-Srivastava2007.tex}, the following necessary and sufficient conditions are obtained: For real numbers $\alpha\ne 0$ and $\beta$, the function
\begin{equation*}
g_{\alpha,\beta}(x)=\biggl[\frac{e^x\Gamma(x+1)} {(x+\beta)^{x+\beta}}\biggr]^\alpha,\quad x\in(\max\{0,-\beta\},\infty)
\end{equation*}
is logarithmically completely monotonic if and only if either $\alpha>0$ and $\beta\geq1$ or $\alpha<0$ and $\beta\leq\frac12$. Further considering the fact in~\cite[p.~98]{Dubourdieu} that a completely monotonic function which is non-identically zero cannot vanish at any point on $(0,\infty)$ gives
$$
(-1)^k[\ln g_{\alpha,\beta}(x)]^{(k)}=(-1)^k\alpha[x+\ln\Gamma(x)+\ln x-(x+\beta)\ln(x+\beta)]^{(k)}>0
$$
for $k\in\mathbb{N}$ and $x\in(0,\infty)$ if and only if either $\alpha>0$ and $\beta\geq1$ or $\alpha<0$ and $\beta\leq\frac12$. As a result, from straightforward calculation and standard arrangement, inequalities \eqref{qi-psi-ineq-beta-1} and~\eqref{qi-psi-ineq-beta-2} follow. Lemma~\ref{qi-psi-ineq-beta-lem} is thus proved.
\end{proof}

\section{Proofs of theorems}

Now we are in a position to prove our theorems.

\begin{proof}[Proof of Theorem~\ref{Open-TJM-2003-thm-1}]
For $x\ne0$, taking the logarithm of $h_{\alpha,y}(x)$ gives
\begin{equation*}
\ln h_{\alpha,y}(x)=\frac{\ln\Gamma(x+y+1)-\ln\Gamma(y+1)}x-\alpha\ln(x+y+1).
\end{equation*}
A direct differentiation yields
\begin{equation}\label{one-partial-der}
\begin{aligned}{}
[\ln h_{\alpha,y}(x)]^{(k)}&= \frac{k!}{x^{k+1}}\sum_{i=0}^k\frac{(-1)^{k-i}x^i\psi^{(i-1)}(x+y+1)}{i!} \\* &\quad-\frac{(-1)^kk!\ln\Gamma(y+1)}{x^{k+1}}-\frac{(-1)^{k-1}(k-1)!\alpha}{(x+y+1)^k}
\end{aligned}
\end{equation}
for $k\in\mathbb{N}$, where $\psi^{(-1)}(x+y+1)$ and $\psi^{(0)}(x+y+1)$ stand for $\ln\Gamma(x+y+1)$ and $\psi(x+y+1)$ respectively. Furthermore, a simple calculation gives
\begin{align*}
\bigl\{x^{k+1}[\ln h_{\alpha,y}(x)]^{(k)}\bigr\}'&= (-1)^{k-1}x^k\biggl[(-1)^{k-1}\psi^{(k)}(x+y+1)\\*
&\quad -\frac{(k-1)!\alpha}{(x+y+1)^{k}}-\frac{k!(y+1)\alpha}{(x+y+1)^{k+1}}\biggr].
\end{align*}
Utilizing~\eqref{qi-psi-ineq} in the above equation leads to
\begin{multline}\label{DB-ineq-h}
\frac{(k-1)!(1-\alpha)}{(x+y+1)^{k}}+\frac{k![1/2-(y+1)\alpha]}{(x+y+1)^{k+1}} \le\frac{(-1)^{k-1}}{x^k}\bigl\{x^{k+1}[\ln h_{\alpha,y}(x)]^{(k)}\bigr\}' \\ \le\frac{(k-1)!(1-\alpha)}{(x+y+1)^{k}}+\frac{k![1-(y+1)\alpha]}{(x+y+1)^{k+1}}
\end{multline}
for $k\in\mathbb{N}$, $x\ne0$, $y\in(-1,\infty)$ and $\alpha\in(-\infty,\infty)$. Therefore,
\begin{equation}\label{2-part-pos-neg}
\frac{(-1)^{k-1}}{x^k}\bigl\{x^{k+1} [\ln h_{\alpha,y}(x)]^{(k)}\bigr\}'
\begin{cases}
\le0,&\text{if $\alpha\ge1$ and $\alpha\ge\frac1{y+1}$}\\
\ge0,&\text{if $\alpha\le1$ and $\alpha\le\frac1{2(y+1)}$}
\end{cases}
\end{equation}
for $k\in\mathbb{N}$, $y>-1$ and $x\ne0$.
For $x>0$, the equation~\eqref{2-part-pos-neg} means
\begin{equation*}
\bigl\{x^{2k} [\ln h_{\alpha,y}(x)]^{(2k-1)}\bigr\}'
\begin{cases}
\le0,&\text{if $\alpha\ge1$ and $\alpha\ge\frac1{y+1}$}\\
\ge0,&\text{if $\alpha\le1$ and $\alpha\le\frac1{2(y+1)}$}
\end{cases}
\end{equation*}
and
\begin{equation*}
\bigl\{x^{2k+1} [\ln h_{\alpha,y}(x)]^{(2k)}\bigr\}'
\begin{cases}
\ge0,&\text{if $\alpha\ge1$ and $\alpha\ge\frac1{y+1}$}\\
\le0,&\text{if $\alpha\le1$ and $\alpha\le\frac1{2(y+1)}$}
\end{cases}
\end{equation*}
for $k\in\mathbb{N}$. From~\eqref{one-partial-der}, it is easy to see that
\begin{equation}\label{0-to-limit}
\lim_{x\to0}\bigl\{x^{k+1}[\ln h_{\alpha,y}(x)]^{(k)}\bigr\}=0
\end{equation}
for $k\in\mathbb{N}$ and any given $y>-1$. As a result,
\begin{equation}\label{h-der->0}
[\ln h_{\alpha,y}(x)]^{(2k-1)}
\begin{cases}
<0,&\text{if $\alpha\ge1$ and $\alpha\ge\frac1{y+1}$}\\>0,&\text{if $\alpha\le1$ and $\alpha\le\frac1{2(y+1)}$}
\end{cases}
\end{equation}
and
\begin{equation}\label{h-der-<0}
[\ln h_{\alpha,y}(x)]^{(2k)}
\begin{cases}
>0,&\text{if $\alpha\ge1$ and $\alpha\ge\frac1{y+1}$}\\
<0,&\text{if $\alpha\le1$ and $\alpha\le\frac1{2(y+1)}$}
\end{cases}
\end{equation}
for $k\in\mathbb{N}$ and $x\in(0,\infty)$, that is,
\begin{equation}\label{h-der-=0}
(-1)^k[\ln h_{\alpha,y}(x)]^{(k)}
\begin{cases}
>0,&\text{if $\alpha\ge1$ and $\alpha\ge\frac1{y+1}$}\\<0,&\text{if $\alpha\le1$ and $\alpha\le\frac1{2(y+1)}$}
\end{cases}
\end{equation}
for $k\in\mathbb{N}$ and $x\in(0,\infty)$. Hence, the function~\eqref{fun} is logarithmically completely monotonic with respect to $x$ on $(0,\infty)$ if $\alpha\ge1$ and $\alpha\ge\frac1{y+1}$, and so is the reciprocal of the function~\eqref{fun} if either $0<\alpha\le1$ and $\alpha\le\frac1{2(y+1)}$ or $\alpha\le0$ and $y>-1$.
\par
If $x\in(-y-1,0)$, the equation~\eqref{2-part-pos-neg} means
\begin{equation*}
\bigl\{x^{k+1} [\ln h_{\alpha,y}(x)]^{(k)}\bigr\}'
\begin{cases}
\ge0,&\text{if $\alpha\ge1$ and $\alpha\ge\frac1{y+1}$}\\
\le0,&\text{if $\alpha\le1$ and $\alpha\le\frac1{2(y+1)}$}
\end{cases}
\end{equation*}
for $k\in\mathbb{N}$. By virtue of~\eqref{0-to-limit}, it follows that
\begin{equation*}
x^{k+1} [\ln h_{\alpha,y}(x)]^{(k)}
\begin{cases}
\le0,&\text{if $\alpha\ge1$ and $\alpha\ge\frac1{y+1}$}\\
\ge0,&\text{if $\alpha\le1$ and $\alpha\le\frac1{2(y+1)}$}
\end{cases}
\end{equation*}
for $k\in\mathbb{N}$, which is equivalent to that the equations~\eqref{h-der->0} and~\eqref{h-der-<0} hold for $x\in(-y-1,0)$. As a result, the equation~\eqref{h-der-=0} is valid for $k\in\mathbb{N}$ and $x\in(-y-1,0)$. Therefore, the function $h_{\alpha,y}(x)$ has the same logarithmically complete monotonicity properties on $(-y-1,0)$ as on $(0,\infty)$.
\par
Conversely, if $h_{\alpha,y}(x)$ is logarithmically completely monotonic on $(-y-1,\infty)$, then $[\ln h_{\alpha,y}(x)]'<0$ on $(-y-1,\infty)$, which can be simplified as
\begin{align}\label{func-19}
\alpha&\ge(x+y+1)\biggl[\frac{1}{x^{2}}\sum_{i=0}^1\frac{(-1)^{1-i}x^i\psi^{(i-1)}(x+y+1)}{i!} +\frac{\ln\Gamma(y+1)}{x^2}\biggr] \\
\begin{split}\label{func-20}
&=\frac1{x^2}[(x+y+1)\ln\Gamma(y+1)-(y+1)(x+y+1)\psi(x+y+1)\\
&\quad+(x+y+1)^2\psi(x+y+1)-(x+y+1)\ln\Gamma(x+y+1)]
\end{split}\\ \label{func-21}
&=\frac{x+y+1}{x}\biggl[\frac{x\psi(x+y+1)-\ln\Gamma(x+y+1)}{x} +\frac{\ln\Gamma(y+1)}{x}\biggr].
\end{align}
From~\eqref{qi-psi-ineq-1}, it is easy to see that
\begin{equation}\label{lim-psi-1}
\lim_{x\to0^+}\bigl[x^2\psi(x)\bigr]=0.
\end{equation}
It is common knowledge that
\begin{equation}\label{gamma-identity-period}
\Gamma(x+1)=x\Gamma(x)
\end{equation}
for $x>0$. Taking the logarithm on both sides of~\eqref{gamma-identity-period}, rearranging, and taking limit lead to
\begin{equation}\label{lim-psi-2}
\lim_{x\to0^+}[x\ln\Gamma(x)]=\lim_{x\to0^+}[x\ln\Gamma(x+1)]-\lim_{x\to0^+}[x\ln x]=0.
\end{equation}
Taking logarithmic derivatives on both sides of~\eqref{gamma-identity-period} yields
\begin{equation}
\psi(x+1)=\frac1x+\psi(x)
\end{equation}
for $x>0$, and so
\begin{equation}\label{lim-psi-3}
\lim_{x\to0^+}[x\psi(x)]=1+\lim_{x\to0^+}[x\psi(x+1)]=1.
\end{equation}
Thus, by utilizing~\eqref{lim-psi-1}, \eqref{lim-psi-2} and~\eqref{lim-psi-3}, it is revealed that the limit of the function~\eqref{func-20} as $x\to(-y-1)^+$, that is, as $x+y+1\to0^+$, equals $\frac1{y+1}$. By L'H\^ospital's rule and the double inequality~\eqref{qi-psi-ineq} for $k=1$, we have
\begin{gather*}
\lim_{x\to\infty}\frac{x\psi(x+y+1)-\ln\Gamma(x+y+1)}{x} =\lim_{x\to\infty}[x\psi'(x+y+1)]=1.
\end{gather*}
Hence, the limit of the function~\eqref{func-21} as $x\to\infty$ equals $1$. In a word, the necessary condition for $h_{\alpha,y}(x)$ to be logarithmically completely monotonic is $\alpha\ge1$ and $\alpha\ge\frac1{y+1}$.
\par
If the reciprocal of $h_{\alpha,y}(x)$ is logarithmically completely monotonic, then the inequality~\eqref{func-19} is reversed. Since the limit of the function~\eqref{func-21} equals~$1$ as $x\to\infty$, as showed above, then the necessary condition $\alpha\le1$ is obtained. The proof of Theorem~\ref{Open-TJM-2003-thm-1} is complete.
\end{proof}

\begin{proof}[Proof of Theorem~\ref{Open-TJM-2003-thm-ineq}]
This follows from the monotonicity properties established in Theorem~\ref{Open-TJM-2003-thm-1}.
\end{proof}

\begin{proof}[Proof of Theorem~\ref{vol-ball-ineq-thm}]
Letting $t=1$, $y=0$ and $x=\frac{n}2$ for $n\in\mathbb{N}$ in~\eqref{open-TJM-2003-ineq} reveals that
\begin{equation*}
\frac{n+2}{n+4}<\frac{[\Gamma(n/2+1)]^{2/n}}
{[\Gamma((n+2)/2+1)]^{2/(n+2)}} <\sqrt{\frac{n+2}{n+4}}\,
\end{equation*}
which is equivalent to the inequality~\eqref{sqrt-frac-n+2n+4}.
\par
If taking $y=1$, $t=1$ and $x=\frac{n}2-1$ for $n\in\mathbb{N}$ in~\eqref{open-TJM-2003-ineq}, then the inequality~\eqref{vol-ball-1} follows.
\par
Replacing $t$ by $\frac12$, $y$ by $0$, and $x$ by $\frac{n}2$ in~\eqref{open-TJM-2003-ineq} and simplifying result in~\eqref{ratio-Omega-n-n+1}.
\end{proof}

\section{Remarks}

After proving our theorems, we give several remarks about them.

\begin{rem}
Theorem~\ref{Open-TJM-2003-thm-1} extends and generalizes the logarithmically complete monotonicity of the function~\eqref{fun-ori} established in~\cite{Extension-TJM-2003.tex} and a part of the results in~\cite{Zhao-Chu-Jiang}.
\end{rem}

\begin{rem}
The inequality~\eqref{open-TJM-2003-ineq} generalizes and extends the inequality~\eqref{new3} and the main results in \cite{Open-TJM-2003-Ineq-Ext-JAT.tex, Ya-Ming-Yu-JMAA-09}: For $x+y>0$ and $y+1>0$ the inequality
\begin{equation}
\frac{[\Gamma(x+y+1)/\Gamma(y+1)]^{1/x}}{[\Gamma(x+y+2)/\Gamma(y+1)]^{1/(x+1)}} <\biggl(\frac{x+y}{x+y+1}\biggr)^{1/2}
\end{equation}
is valid if $x>1$ and reversed if $x<1$ and that the power $\frac12$ is the best possible.
\end{rem}

\begin{rem}
When $n>2$, the inequality~\eqref{ratio-Omega-n-n+1} refines the following double inequality in~\cite[Theorem~1]{ball-volume-rn}:
\begin{equation}\label{Theorem-1-ball-volume-rn}
\frac2{\sqrt{\pi}\,}\Omega_{n+1}^{n/(n+1)}\le\Omega_n <\sqrt{e}\,\Omega_{n+1}^{n/(n+1)}, \quad n\in\mathbb{N}.
\end{equation}
For more information on inequalities for the volume of the unit ball $\mathbb{B}^n$ in $\mathbb{R}^n$, please refer to \cite{alzer-ball-ii, Berg-Pedersen-ball-PAMS, Berg-Pedersen-ball-ArXiv, unit-ball.tex}, Section~7.5 in~\cite[pp.~72\nobreakdash--73]{bounds-two-gammas.tex} and related references therein.
\end{rem}

\begin{rem}
Theorem~\ref{Open-TJM-2003-thm-1} may be restated as follows: For $y\in(0,\infty)$ and $\alpha\in(-\infty,\infty)$, the function
\begin{equation}\label{fun-y+1}
H_y(x)=
\begin{cases}
\dfrac{[{\Gamma(x+y)}/{\Gamma(y)}]^{1/x}}{(x+y)^\alpha}, &x\in(-y,\infty)\setminus\{0\}\\[0.8em]
\dfrac{e^{\psi(y)}}{y^\alpha},&x=0
\end{cases}
\end{equation}
is logarithmically completely monotonic with respect to $x\in(-y,\infty)$ if and only if $\alpha\ge\max\bigl\{1,\frac1{y}\bigr\}$, and so is its reciprocal if $\alpha\le\min\bigl\{1,\frac1{2y}\bigr\}$ and only if $\alpha\le1$.
\end{rem}

\begin{rem}
We conjecture that when $y>-\frac12$ the condition $\alpha\le\frac1{2(y+1)}$ is also necessary for the reciprocal of the function~\eqref{fun} to be logarithmically completely monotonic with respect to $x\in(-y-1,\infty)$.
\end{rem}

\begin{rem}
For more information on the history, backgrounds, motivations and recent developments of the topic in this paper, please refer to \cite{Abram-Baric-Mat-Pecaric-JIA-07, Bennett-JIA-07, cqcd, Guo-Qi-TJM-03.tex, rgmg-2.tex, Sqrt[n]{Frac(n+k)!{k!}}.tex-soochow, Extension-TJM-2003.tex, msicfs.tex-mia, Martin-JIA-07, Minc-Sathre.tex-tamkang, Qi-Sun-Anal-Math-2006} and related references therein.
\end{rem}

\begin{rem}
In passing, we survey the history of the notion ``logarithmically completely monotonic function''. By searching for the term ``logarithmically completely monotonic function'' in the database \href{http://www.ams.org/mathscinet/}{MathSciNet}, it is found that this phrase was probably first used in~\cite{Atanassov}, but without an explicit definition. Thereafter, it seems to have not been used by the mathematical community. In early 2004, this terminology was again used in~\cite{minus-one} (the preprint of \cite{compmon2, minus-one-JKMS.tex}) and it was immediately referenced in~\cite{grin-ismail} and \cite{auscm-rgmia} (the preprint of~\cite{e-gam-rat-comp-mon}). In \cite[Theorem~4]{minus-one}, it was proved that a logarithmically completely monotonic function $f(x)$ on $I$ must be completely monotonic (i.e., the inequality
$
0\le(-1)^{k}f^{(k)}(x)<\infty
$
holds for all $k\geq0$ on $I$), but not conversely. This result was announced while revising~\cite{compmon2}. This conclusion and its proofs were presented once and again in~\cite{CBerg} and~\cite{schur-complete} (the preprint of~\cite{absolute-mon-simp.tex}). More importantly, in the paper~\cite{CBerg}, the logarithmically completely monotonic functions on $(0,\infty)$ were characterized as the infinitely divisible completely monotonic functions studied in~\cite{horn} and all Stieltjes transforms were proved to be logarithmically completely monotonic on $(0,\infty)$. For information on the completely monotonic functions, please refer to~\cite[Chapter~XIII]{mpf-1993} and~\cite[Chapter~IV]{widder}, especially to the recently published monograph~\cite{Schilling-Song-Vondracek-2010}.
\end{rem}

\begin{rem}
This paper is a part of the preprint~\cite{Open-TJM-2003.tex}.
\end{rem}

\end{document}